\documentclass[12pt,psfig]{article}

\usepackage{amsfonts,amsmath,leqno,graphicx}
\usepackage{epsfig,subfigure,fancybox}
\usepackage{amsmath}
\usepackage{amssymb}
\usepackage{amsfonts}
\usepackage{graphicx}
\usepackage{dsfont}
\usepackage[subfigure]{graphfig}

\usepackage{slashbox}

\newcommand{\rr}{{\mathbb R}}

\renewcommand{\thesection}{\arabic{section}}
\renewcommand{\theequation}{\thesection.\arabic{equation}}
\newcommand{\beq}[1]{\begin{equation} \label{#1}}
\newcommand{\eeq}{\end{equation}}
\newcommand{\bea}{\bed\begin{array}{rl}}
\newcommand{\eea}{\end{array}\eed}
\newcommand{\bed}{\begin{displaymath}}
\newcommand{\eed}{\end{displaymath}}
\newcommand{\barray}{\begin{array}{ll}}
\newcommand{\earray}{\end{array}}

\newcommand{\disp}{\displaystyle}

\newcommand{\al}{\alpha}

\newcommand{\NM}{{\rm NM}}
\newcommand{\AM}{{\rm AM}}

\newcommand{\lbar}{\overline}
\newcommand{\wdt}{\widetilde}
\newcommand{\wdh}{\widehat}

\def\para#1{\vskip 0.4\baselineskip\noindent{\bf #1}}

\def\rr{{\Bbb R}}
\def\({\left(}
\def\){\right)}
\def\one{{\hbox{1{\kern -0.35em}1}}}

\begin{document}

\title{Deep Filtering\thanks{This research  was supported in part by the Army Research Office under grant W911NF-19-1-0176.} }

\author{Q. Zhang,\thanks{Department of Mathematics,
University of Georgia, Athens, GA 30602, qz@uga.edu.}\and
G. Yin,\thanks{Department of Mathematics, Wayne
  State University, Detroit, MI 48202, gyin@wayne.edu.}\and
L.Y. Wang\thanks{Department of Electrical and Computer Engineering,
  Wayne State University, Detroit, MI 48202, {\mbox lywang@wayne.edu}.}}


\date{ }

\maketitle

\begin{abstract}
  This paper develops a deep learning method for linear and nonlinear
  filtering. The idea is to start with a nominal dynamic model and
  generate Monte Carlo sample paths. Then these samples are used to train
  a deep neutral network. A least square error is used as a loss function for
  network training. Then the resulting weights are applied to Monte Carlo samples
  from an actual dynamic model. The deep filter obtained in such a way compares
  favorably to the traditional Kalman filter in linear cases and
  the extended Kalman filter in nonlinear cases. Moreover,
   a switching model with jumps is studied to show the adaptiveness and power
of our deep filtering method.
  A main advantage of deep filtering is its robustness
  when the nominal model and actual model differ.
  Another advantage of deep filtering is that real data can be used directly to train
  the deep neutral network. Therefore,
   one does not need to calibrate the  model.


\end{abstract}

\maketitle

\section{Introduction}
This paper develops  deep learning methods for both linear and nonlinear filtering problems. Recent advent of applications of artificial intelligence in diversified domains has promoted
an intensified interest in using machine learning theory to
treat stochastic dynamic systems and stochastic controls.
It opens up many possibilities in state estimation
with reduced computational complexity, alleviating
the curse of dimensionality.
There are numerous successful applications of
deep  learning in multi-agent systems,  traffic control,
robotics, personalized recommendations, and games of GO and Atari.
Despite many progresses, there seems to be no work devoted to using the deep learning approach in
state estimation and filtering.
This paper aims to develop deep neutral network based
filtering schemes.

\begin{center}
\begin{figure}[htb]
\begin{center}
\mbox{
         \psfig{figure=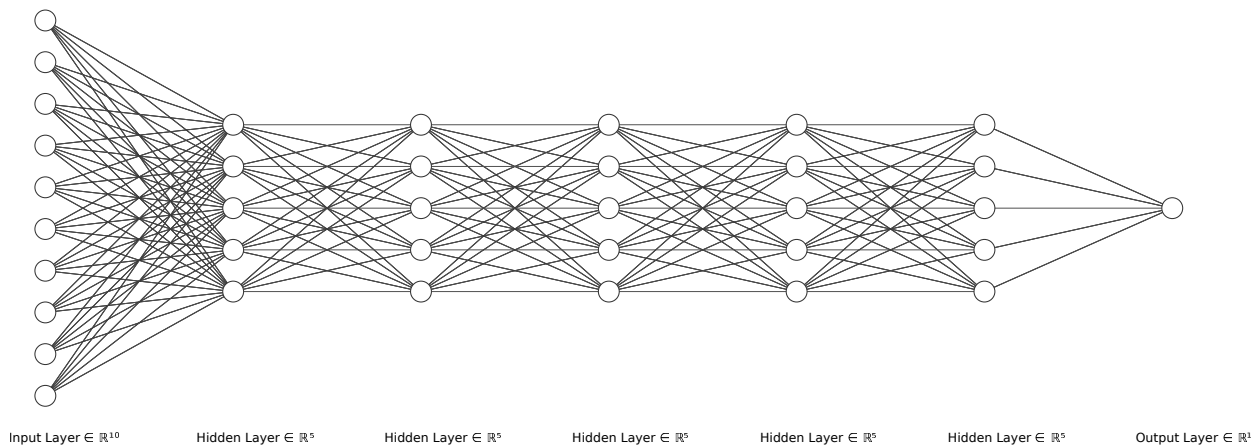,width=1.0\linewidth }}
    \end{center}
\caption{A Deep Neutral Network}\label{DNN}
\end{figure}
    \end{center}

\para{Deep Neutral Networks (DNN) and Backpropagation.}
Neural networks (NN) are often used to approximate functions, which
are often complex and highly nonlinear, arising from a wide variety of applications.
The main approaches
are of compositional nature
and rely on composition of hidden layers of base functions.
A deep neural network is in fact,
an NN with several hidden layers. The deepness of the network is measured by the number of layers used.
In this paper, we only consider a fully connected NN and there are no
connections between nodes in the same layer.
A typical DNN diagram of such a class is depicted in Figure~\ref{DNN}.
For related literature on DNN, we refer the reader to the online book
by Nielsen \cite{Nielsen}.

The backpropagation is a commonly used main driver in DNN training.
By and large, backpropagation is an algorithm for supervised learning
of artificial neural networks using  gradient descent procedures.
Given an artificial NN and a loss or error function, the scheme
calculates the gradient of the loss function
with respect to the weights of the neural network.
The
calculation of the gradient
proceeds backwards through the network starting with the gradient of the final layer of weights.
To facilitate the computation,
partial computations of the gradient from one layer
are reused for the previous layer's gradient calculation.
Such backwards flow of information is designed for
efficient computation of the gradient at each layer.
In particular, backpropagation requires three items:
\begin{itemize}
\item[(a)] A data set consisting of fixed pairs of input and output variables;
\item[(b)] A feedforward NN with parameters given by the weights $\theta$;
\item[(c)] A loss (error) function $L(\theta)$ defining the error
between the desired output and the calculated output.
\end{itemize}
In this paper, the NN training will use the stochastic gradient decent method to find
the weight vector $\theta$ to minimize a loss function $L(\theta)$.
The details are to be given later.

\para{Linear and Nonlinear Filtering.}
As is well known,
filtering is concerned with dynamic systems in which the internal state variables
are not completely observable. There are numerous real-world applications
involving state estimation and
filtering, including maneuvered target tracking, speech recognition,
telecommunications, financial engineering, among many others.
Traditional approaches derive estimators based on observations
with least square errors.
Working at a setup in discrete time for the underlying systems and
under a Gaussian distribution framework,
the corresponding filtering problem is to find the
conditional mean of the state given the observation up to time $n$.
The best known filter is the Kalman filter for linear models.
As for
 some nonlinear models,  extended Kalman filters can be applied.
We refer the reader to Anderson and Moore \cite{AndersonM} for
more details.

Early development in nonlinear filtering can be found in
Duncan \cite{Duncan}, which focuses on conditional densities for diffusion
processes, Mortensen \cite{Mortensen}, which concentrates on the
most probable trajectory approach;
Kushner \cite{Kushner}, which derives nonlinear filtering equations,
and Zakai \cite{Zakai}, which uses unnormalized equations.

There are many progresses made in the past decades since then.
For example, hybrid filtering can be found in
Hijab \cite{Hijab-paper} with an unknown constant,
Zhang \cite{Zhang-switch} with small observation noise,
Miller and Runggaldier \cite{MillerR} with Markovian
jump times,
Blom and Bar-Shalom \cite{BlomB} for discrete-time hybrid model
and the Interactive Multiple Model algorithm,
Dufour et al \cite{DufourBE1,DufourBE2} and
Dufour and Elliott \cite{DufourE} for models with regime switching.
Some later developments along this line can be found in
Zhang \cite{Zhang-switch,Zhang1,Zhang2}.
Despite these important progresses, the computation of filtering remains a daunting
task. For nonlinear filtering, there have been yet feasible and efficient  schemes to mitigate
high computational complexity (with infinite dimensionality). Much effort has been devoted to finding
 computable approximation schemes.

\para{Deep Filtering.}
In this paper, we propose a new framework, which
focuses on deep neutral network based filtering.
Under a given model, the idea is to generate Monte Carlo samples
and then use these samples to train a deep neutral network.
The observation process generates
inputs to the DNN;
the state
from the Monte Carlo samples is used as the target.
A least square error of the target and calculated output
is used as a loss function for network training
for
weight vectors. Then these weight vectors
are applied to
another set of Monte Carlo samples of the actual dynamic model.
The corresponding calculated DNN output is
termed
a deep filter (DF).

In this paper, we demonstrate the adaptiveness, robustness, and effectiveness of
our DF through numerical examples.
The deep filter compares favorably
to the traditional Kalman filter in linear cases and
  the extended Kalman filter in nonlinear cases.
  Moreover,
   a switching model with jumps is studied and  used to show the
  feasibility and
  flexibility of our deep filtering.
  A major advantage of deep filter is its robustness
  when the nominal model and actual model differ.
  Another advantage of the deep filtering is that
  real data can be used directly to train
  the deep neutral network. Therefore, model calibration is no longer needed
  in applications.

The rest of the paper is arranged as follows.
Section \ref{sec:DF} begins with our deep filtering algorithm,
followed by its corresponding versions for  linear models, nonlinear models, and
switching models.
Numerical experiments are presented.
Performance of the deep filter is examined through various models
and compared with benchmark linear and nonlinear filters.
Concluding remarks are provided in Section \ref{sec:con}.

\section{
Deep Filter}\label{sec:DF}

Let $x_n \in \rr^{m_1}$ denote the state process with system equation
\[
x_{n+1}=F_n(x_n,u_n),\ x_0=x,\ n=0,1,2,\ldots,
\]
for some suitable functions $F_n: \rr^{m_1} \times \rr^{l_1} \mapsto \rr^{m_1} $ and system noise $\{u_n\}$ with $u_n \in \rr^{l_1}$.
A function of $x_n$ can be observed with possible noise corruption.
In particular, the observation process $y_n\in \rr^{m_2}$ is given by
\[
y_n=H_n(x_n,v_n),
\]
with noise $\{v_n\}$, $v_n \in \rr^{l_2}$, and $H_n: \rr^{m_1} \times \rr^{l_2} \mapsto \rr^{m_2}$.   Next, we propose our
{\bf deep filter} as follows.

Let $N_{\rm seed}$ denote the number of training sample paths and
let $n_0$ denote the training window for each sample path.
For any fixed $\kappa=n_0,\ldots,N$
with a fixed
$\omega$, we take
$\{y_\kappa(\omega),y_{\kappa-1}(\omega),\ldots,y_{\kappa-n_0+1}(\omega)\}$
as the input vector to the neural network and $x_\kappa(\omega)$ as the target. In what follows, we shall suppress the $\omega$ dependence.
Fix $x_\kappa$,
let
$\xi_\ell$
denote the neural network output at iteration $\ell$, which depends on the parameter $\theta$ as well as a noise $\{\zeta_\ell\}$.
The noise $\{\zeta_\ell\}$ collects all the random influence in the filtering process. A simplest form of $\xi_\ell(\theta,\zeta_\ell) = \xi(\theta)+ \zeta_\ell$ (i.e., $\xi_\ell$ is independent of $\ell$ and the noise is additive). The formulation here, however, includes more general cases as possibilities.
Our goal is to
find an NN weight $\theta \in \rr^{m_3}$ so as to minimize
the loss function
\beq{loss}
L(\theta)={1\over 2} E
 |\xi_\ell(\theta,\zeta_\ell)- x_\kappa|^2.
\eeq
Recall that we do this for fixed $x_\kappa$.
We follow the backpropagation method to search the optimal weights.
The stochastic gradient decent will be used throughout, which takes the form
\beq{sa-N}
\theta_{\ell+1}=\theta_\ell-\gamma {\partial \xi_\ell(\theta_\ell,\zeta_\ell) \over \partial \theta} [\xi_\ell (\theta_\ell, \zeta_\ell) -x_\kappa ]  ,\eeq
with learning rate $0<\gamma<1$.
Note that  $(\partial /\partial \theta) \xi_\ell(\theta,\zeta)$ is a matrix of the dimension ${m_3\times m_1}$.
Define $\theta^\gamma(t)= \theta_\ell$ for $t \in [\ell \gamma, \ell\gamma + \gamma)$.
Under our neural network, it is easy to see the continuously dependent of
the output $\xi_\ell$ on the weight vector $\theta$.
Assume the following average condition: For each $\theta$ and each positive integer $m$,
\[
\frac{1}{\ell}\sum_{j=m}^{\ell+m-1} E_m\xi_\ell(\theta,\zeta_j) \to \bar\xi (\theta) \ \hbox{ in probability as }\ell \to \infty ,
\]
where $E_m$ denotes the conditional expectation on the information up to $m$. Note that we only need a weak law of large number type condition holds as above.
Then it can be shown that $\theta^\gamma(\cdot)$
converges weakly to $\theta(\cdot)$ such that $\theta(\cdot)$ satisfies the differential equation
\beq{ode-lim}
\dot \theta(t)= -  {\partial \bar\xi (\theta (t)) \over \partial \theta} [\bar \xi(\theta(t))- x_\kappa].\eeq
Assume also that there is a $\theta^*$ satisfying
\beq{full-r} [(\partial /\partial \theta) \bar \xi(\theta^*)]'[(\partial /\partial \theta) \bar \xi(\theta^*)] \ \hbox{ is of full rank} \eeq  (i.e., the matrix is an $m_1\times m_1$ full rank matrix).
Then
$${\partial \bar\xi (\theta^* ) \over \partial \theta} [\bar \xi (\theta^*)- x_\kappa] =0,$$
leads to
$$\left({\partial \bar\xi (\theta^* ) \over \partial \theta}\right)'
\left({\partial \bar\xi (\theta^* ) \over \partial \theta}\right) [\bar \xi (\theta^*)- x_\kappa]=0,$$
where $A'$ denotes the transpose of $A$.  
Using  \eqref{full-r} and multiplying  the above equation by
$$\left\{\left({\partial \bar\xi (\theta^* ) \over \partial \theta}\right)'
\left({\partial \bar\xi (\theta^* ) \over \partial \theta}\right)\right\}^{-1}$$
leads to that 
the stationary point $\theta^*$  is given by
$\bar \xi(\theta^*)- x_\kappa=0$. That is the parameter we are searching for is a root of the equation
 $\bar \xi(\theta^*)=x_\kappa$.
Under additional conditions, we can further show that
$\theta^\gamma(\cdot + t_\gamma)$ converges weakly to $\theta^*$ as $\gamma\to 0$, where $t_\gamma\to \infty$ as $\gamma\to 0$.
Since our main effort is to present the deep filtering results, we will not touch upon the convergence of the stochastic gradient algorithm in this paper.

Then these weights $\theta$ are used to out-of-sample
data
$\{\check\omega\}$ with
the actual observation $y_n(\check\omega)$ as inputs in
the subsequent testing stage which leads to neural network output $\wdt x_n(\check\omega)$.
In this paper, $\{\wdt x_n\}$ is called the deep filter.

Note that the training stage is the most time-consuming part.
Normally, it takes a few thousands samples to train the network.
The good part is that such computationally heavy stage is done off-line.
The feedforward part is simple and fast.

 Throughout the rest of this paper, we use numerical examples to evaluate the performance of the deep filter under
various models. We compare its performance with the Kalman filter in linear models and
the extended Kalman filter in nonlinear cases. We also study more general switching
models with jumps and demonstrate the adaptiveness and effectiveness
of the deep filter.

\subsection{Linear Systems}
This section is devoted to
 linear systems.
Let $x_n$ be an $m_1$-dimensional state vector and $y_n$
an $m_2$-dimensional observation vector satisfying the equations:
\beq{linear-sys}
\left\{\begin{array}{l}
x_{n+1}=F_n x_n+G_n u_n,\ x_0=x,\\
y_n=H_n' x_n+v_n, \ n=0,1,2,\ldots,
\end{array}\right.
\eeq
for some matrices $F_n$, $G_n$, and $H_n$ of appropriate dimensions.
Here $u_n$ and $v_n$ are independent random vectors that have Gaussian distributions 
with mean zero and
$E(u_nu_l')=Q^0_n\delta_{nl}$, $E(v_nv_l')=R^0_n\delta_{nl}$,
for $n,l=0,1,2,\ldots$, where $\delta_{nl}=1$ if $n=l$ and $0$ otherwise.

Let ${\mathcal Y}_n=\sigma\{y_j:\ j\leq n\}$
be the  filtration generated by observations and
$\wdh x_n=E[x_n|{\mathcal Y}_{n-1}]$ be the conditional mean.
The corresponding Kalman filter (see Anderson and Moore \cite{AndersonM})
is given by
\beq{KF}
\left\{\begin{array}{l}
\wdh x_{n+1}=F_n\wdh x_n+K_n(y_n-H_n'\wdh x_n),\ \wdh x_0=Ex_0,\\
K_n=F_nR_nH_n(H_n'R_nH_n+R^0_n)^{-1},\\
R_{n+1}=F_n[R_n-R_nH_n(H_n'R_nH_n+R^0_n)^{-1}H_n'R_n]F_n'+G_nQ^0_nG_n',\\
R_0=E[(x_0-\widetilde x_0)(x_0-\widetilde x_0)'].
\end{array}\right.
\eeq
The conditional expectation of $x_n$ given ${\mathcal Y}_n$ can be evaluated
in terms of $\wdh x_n$ 
and $R_n$ as follows:
\[
\lbar x_n=E[x_n|{\mathcal Y}_n]=
\wdh x_n+R_nH_n(H_n'R_nH_n+R^0_n)^{-1}(y_n-H_n'\wdh x_n).
\]

\subsubsection*{Dependence on the NN  Hyperparameters}
We consider the following one-dimensional system:
\beq{general-model}
\left\{\begin{array}{l}
x_{n+1}=(1+0.1\eta)x_n+\sqrt{\eta}\ \sigma u_n,\ x_0=x,\\
y_n= x_n+\sigma_0 v_n, \ n=0,1,2,\ldots,
\end{array}\right.
\eeq
with $u_n$ and $v_n$ being independent Gaussian $N(0,1)$ random variables.

Using window size $n_0=50$ (number of units of input layer) and
$N_{\rm seed}=5000$ to train the DNN as shown in Figure~\ref{DNN}.
The network has 5 hidden layers and each layer has 5 units (neurons).
It has a single output layer.
Also, for all hidden layers, we use the
sigmoid
activation function $\phi(x)=1/(1+e^{-x})$ and  the simple
activation $\phi(x)=x$ for the output layer.
We use the stochastic gradient decent algorithm
with learning rate $\gamma=0.1$.

We select time horizon $T=5$ and step size $\eta=0.005$.
Therefore, the total number of steps is $N=1000$.
We also take $\sigma=0.7$, $\sigma_0=0.5$, and $x_0=1$.
We keep these specifications of the parameters
in the rest of this paper unless otherwise stated.

\begin{center}
\begin{figure}[h!tb]
\begin{center}
\mbox{\subfigure[State $x_n$ (sample path 1)]{
         \psfig{figure=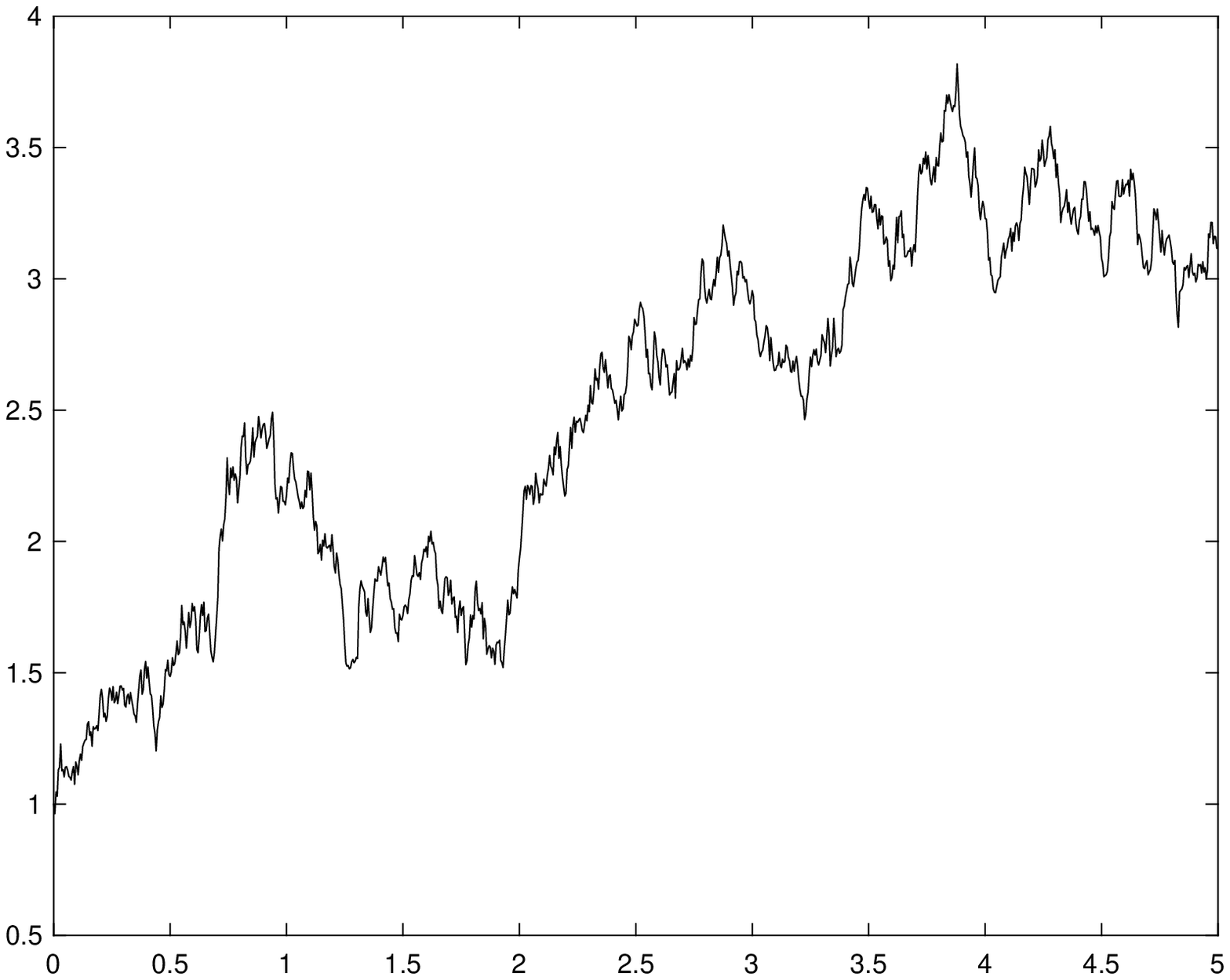,width=0.35\linewidth }} \qquad
      \subfigure[State $x_n$ (sample path 2)]{
        \psfig{figure=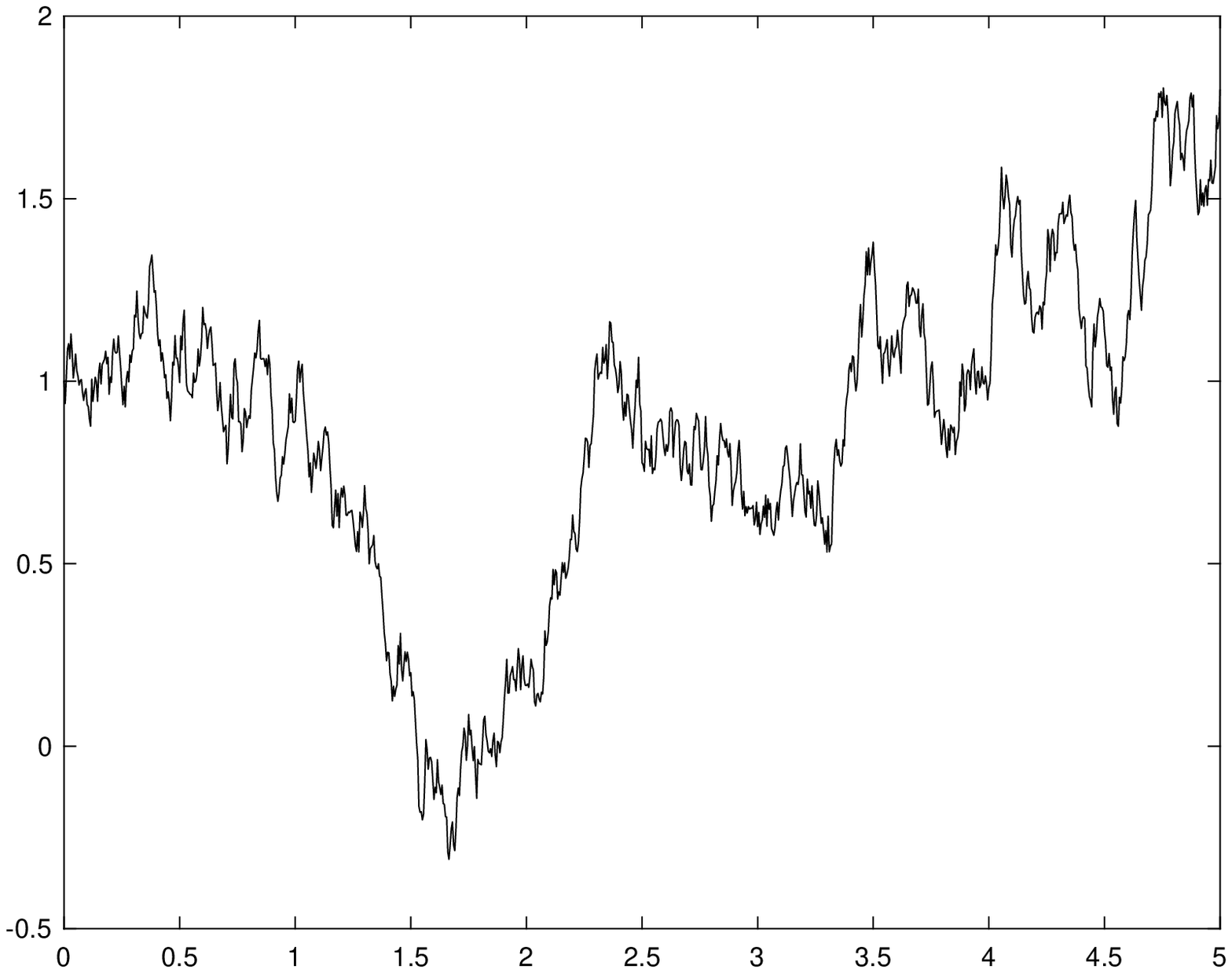,width=0.35\linewidth }}}

\mbox{\subfigure[KF $\wdh x_n$ (sample path 1)]{
         \psfig{figure=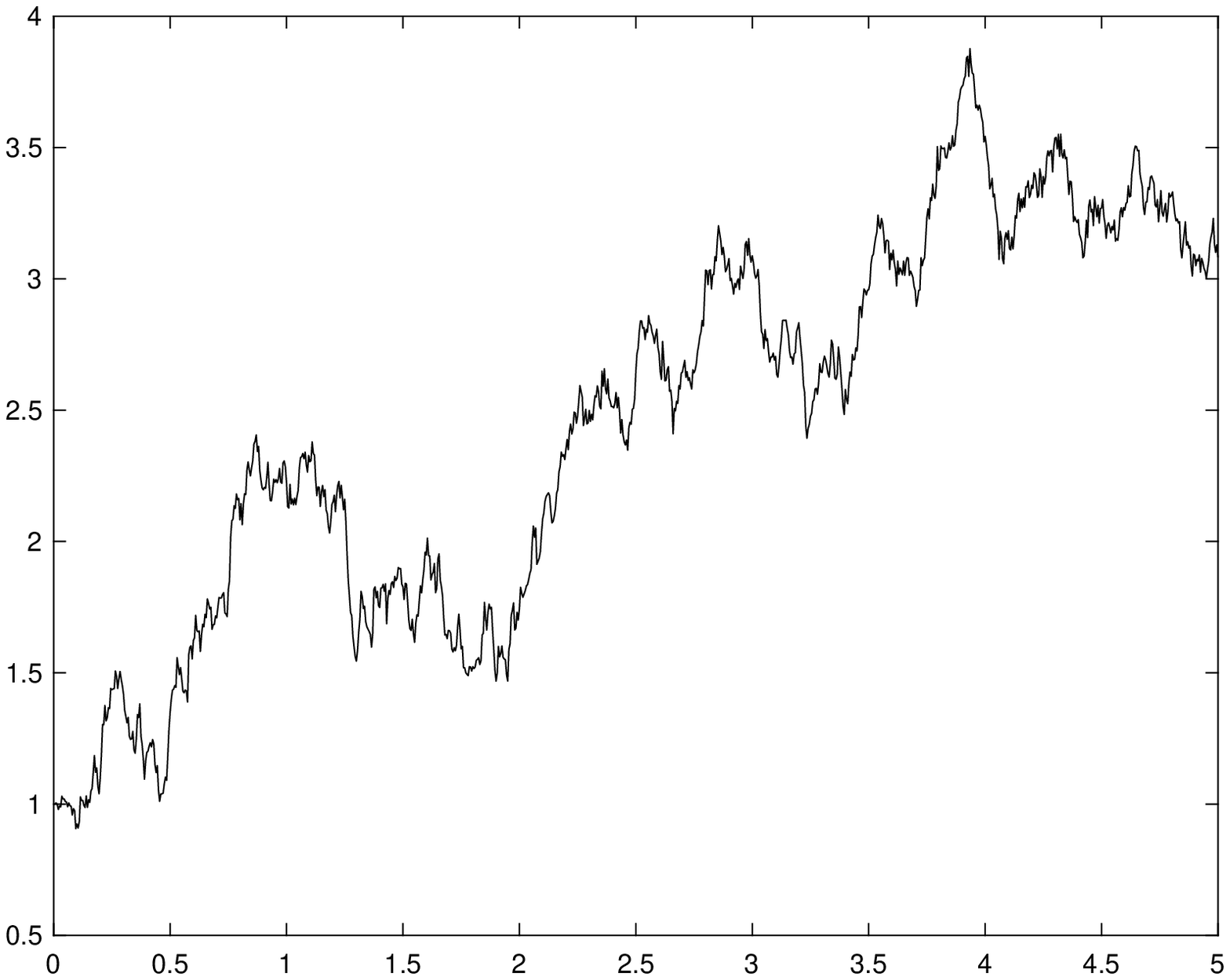,width=0.35\linewidth }} \qquad
      \subfigure[KF $\wdh x_n$ (sample path 2)]{
        \psfig{figure=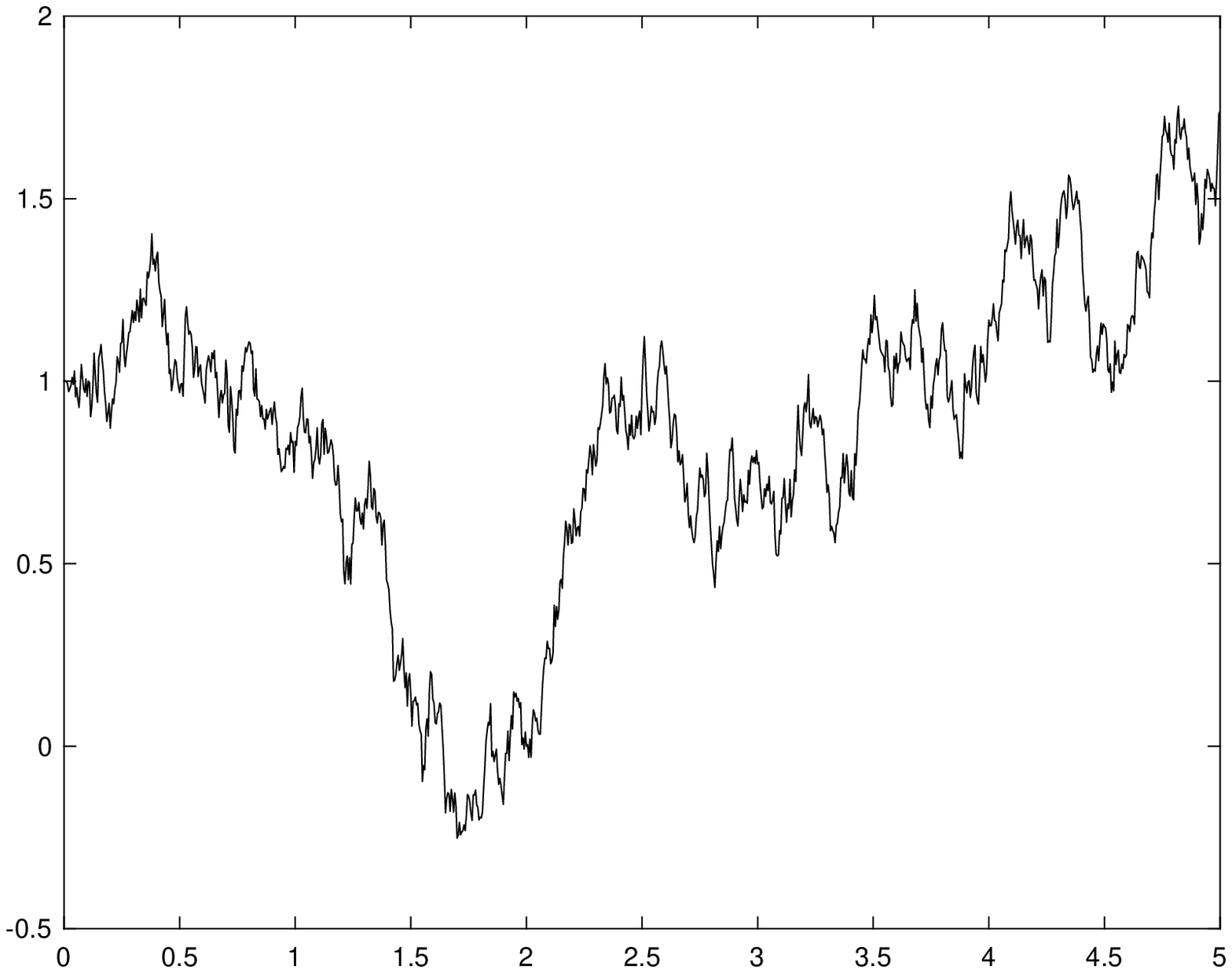,width=0.35\linewidth }}}

\mbox{\subfigure[DF $\wdt x_n$ (sample path 1)]{
         \psfig{figure=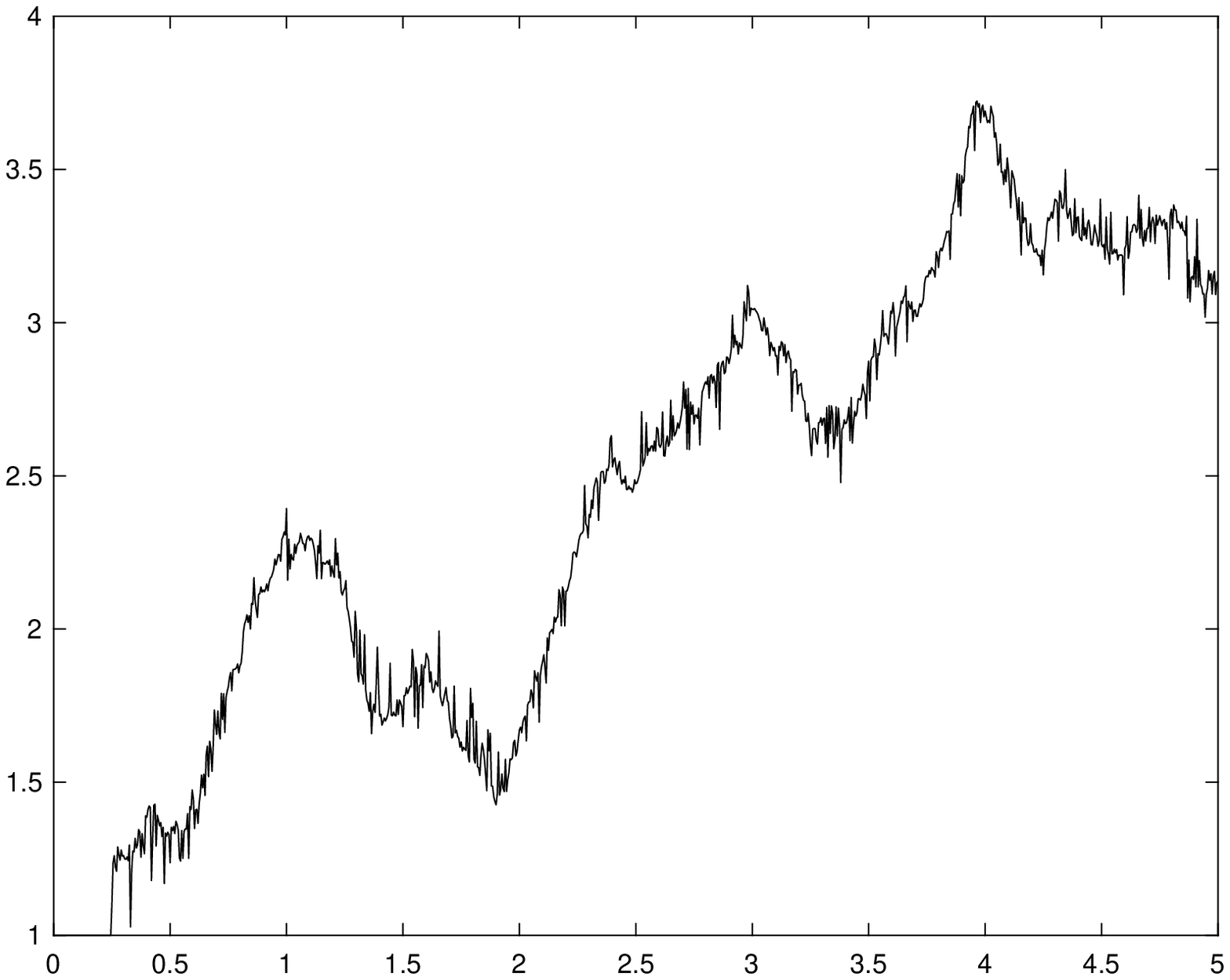,width=0.35\linewidth }} \qquad
      \subfigure[DF $\wdt x_n$ (sample path 2)]{
        \psfig{figure=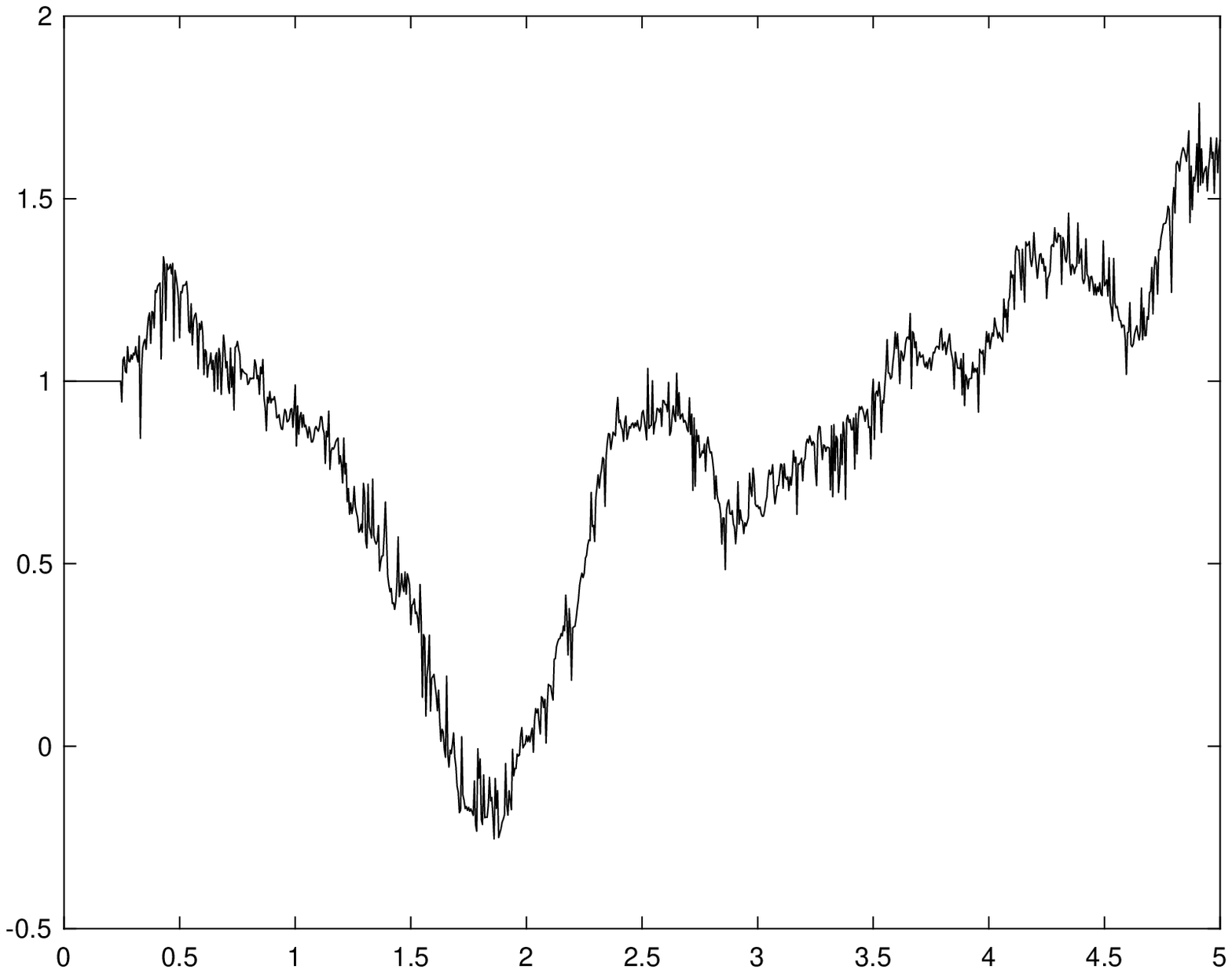,width=0.35\linewidth }}}

\mbox{\subfigure[Error $x_n-\wdt x_n$ (sample path 1)]{
         \psfig{figure=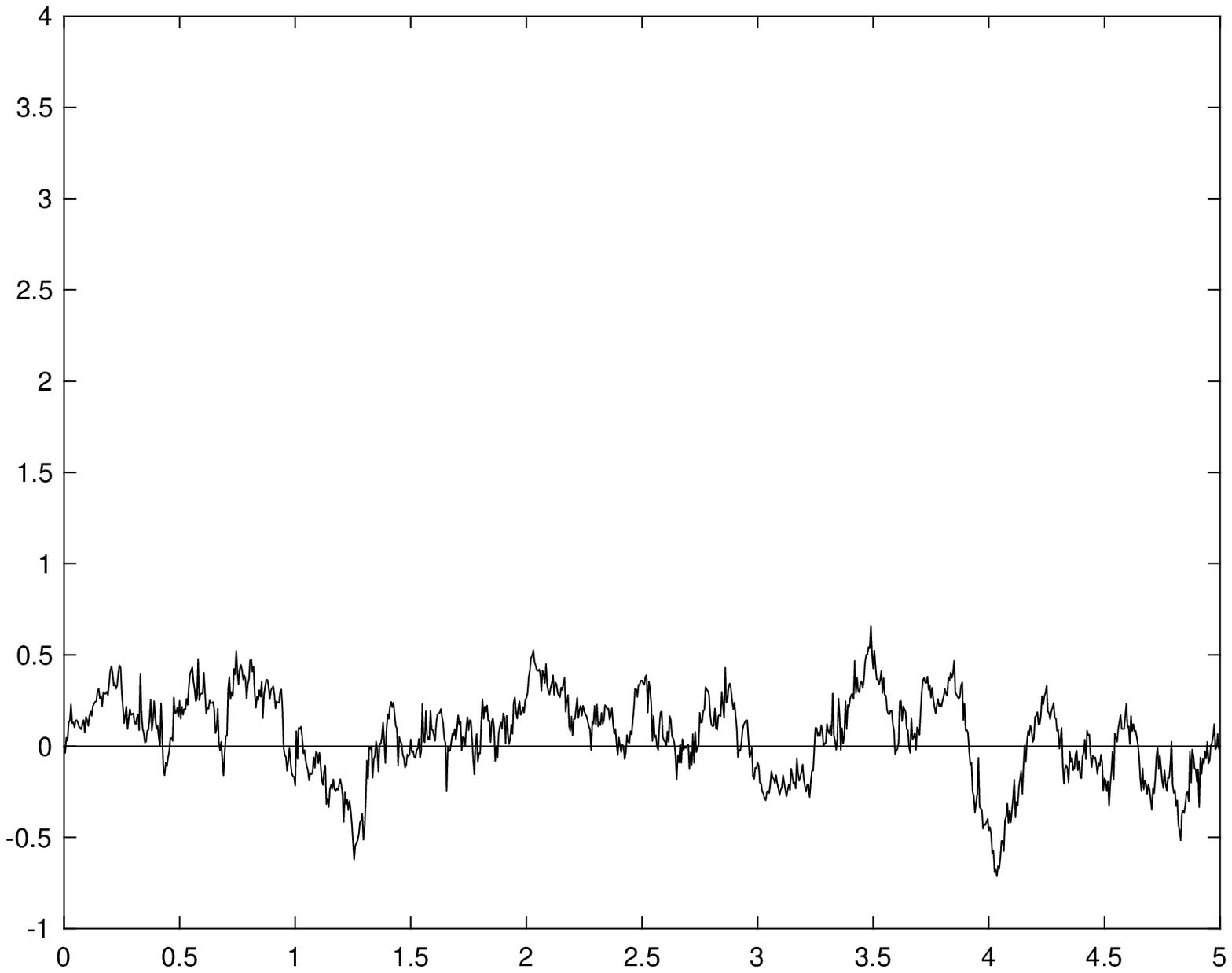,width=0.35\linewidth }} \qquad
      \subfigure[Error $x_n-\wdt x_n$ (sample path 2)]{
         \psfig{figure=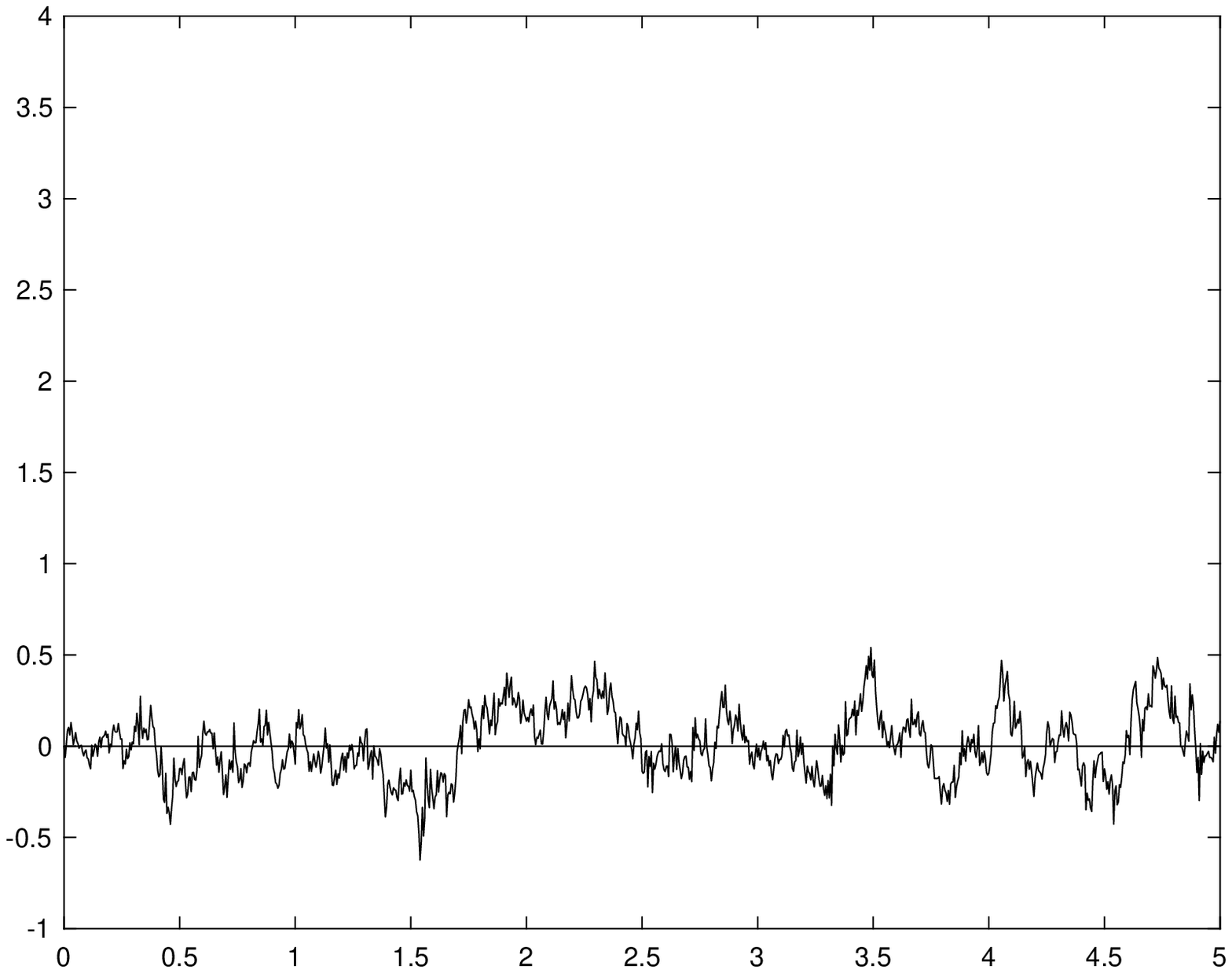,width=0.35\linewidth }}}
\end{center}
\medskip
\caption{\footnotesize Kalman filter vs. Deep filter for linear system:
Two sample paths of $x_n$, $\wdh x_n$, $\wdt x_n$, and errors.}\label{linear-sys}
  \end{figure}
\end{center}

The relative error of vectors
$\xi^1(\check\omega)=(\xi^1_{n_0}(\check\omega),\ldots,\xi^1_N(\check\omega))$
and $\xi^2(\check\omega)=(\xi^2_{n_0}(\check\omega),\ldots,\xi^2_N(\check\omega))$
is defined as
\[
\Vert \xi^1-\xi^2\Vert=
\frac{\disp \sum_{n=n_0}^{N}\sum_{m=1}^{N_{\rm seed}}
  |\xi^1_n(\omega_m')- \xi^2_n(\omega_m')|}{N_{\rm seed}(N-n_0+1){\Sigma}},
\]
where
\[
{\Sigma}=\frac{\disp\sum_{n=n_0}^{N}\sum_{m=1}^{N_{\rm seed}}
(|\xi^1_n(\omega_m')|+|\xi^2_n(\omega_m')|)}{N_{\rm seed}(N-n_0+1)}.
\]
Under these measurements,
we obtain the KF relative error to be
$4.09\%$ and the DF relative error
to be $5.71\%$.
Two sample paths of $x_n$, $\lbar x_n$, $\wdt x_n$, and the corresponding errors are
plotted in Figure~\ref{linear-sys}.

Next, we observe that the
change in the DNN number of layers does not affect much the approximation errors.
For example, when the number of hidden layers changed from 5 to 20, the
corresponding DF relative error changed from $5.71\%$ to $5.76\%$.

Now, we fix the number of hidden layers and vary the number of units of each layer.
The dependence of the DF errors on the number of units for each hidden
 layer is given in Table~\ref{unit-linear}.
    In this table,
    when the observation noise $\sigma_0$ is small,
    a larger number of units
    leads to better approximation. On the other hand, when the observation noise is large,
    this is reversed, i.e., a larger number of units turns out to
    raise the level of the DF.
    Kalman filter errors are also included in this table
    for comparison.
    \begin{table}
    \begin{center}
    {\small
      \begin{tabular}{|c|c|c|c|c|c|}\hline
  \backslashbox{$\hspace*{0.15in}\sigma_0$}{\# NN units}
  &\makebox[3em]{3} &\makebox[3em]{5}&\makebox[3em]{10}
  &\makebox[3em]{20} & \makebox[5em]{KF error}\\ \hline
  0.1 & 2.60  & 2.44 & 2.27  & 2.10  & 1.70  \\ \hline
  0.5 & 5.69  & 5.71 & 4.88  & 4.80  & 4.08  \\ \hline
  1.0 & 6.63  & 6.68 & 6.69  & 7.01  & 5.81  \\ \hline
  2.0 & 9.33  & 9.46 & 9.59  & 10.59  & 8.23  \\ \hline
\end{tabular}
    }
    \end{center}
\smallskip
\caption{Dependence of DF relative errors on \# units of each hidden layer.}\label{unit-linear}
    \end{table}

In Table~\ref{cpu-linear}, the corresponding CPU times (in seconds) are
    provided. Overall, as the number of units increases, the required
    computational time increases. In addition, there appears to be a
    small decrease in CPU time as observation noise increases.
\begin{table}
\begin{center}
    {\small
  \begin{tabular}{|c|c|c|c|c|}\hline
  \backslashbox{$\hspace*{0.15in}\sigma_0$}{\# NN units}
  &\makebox[3em]{3}&\makebox[3em]{5}&\makebox[3em]{10}
  &\makebox[3em]{20} \\ \hline
  0.1 & 31 & 45& 85 & 209   \\ \hline
  0.5 & 30 & 45& 83 & 197  \\ \hline
  1.0 & 28 & 40& 79 & 188   \\ \hline
  2.0 & 27 & 39& 74 & 182   \\ \hline
\end{tabular}
    }
\medskip
    \caption{CPU time (seconds) for DF.}\label{cpu-linear}
    \end{center}
\end{table}

\subsubsection*{Robustness of Deep Filtering}
In this section, we examine the robustness of deep filtering.
We consider separately the nominal model and the actual model.
A nominal model (NM) is an estimated model.
It deviates from real data for different applications.
In this paper, it is used to train our DNNs, i.e.,
a selected mathematical model is used to
generate Monte Carlo sample paths to train the DNN.
The coefficients of the mathematical model are also used
in Kalman filtering equations for comparison.

In real world applications, the conversion from real data to mathematical
models then Monte Carlo processes can be skipped.
Namely, a nominal model consists of actual data to be used to
train the DNN directly.

An actual model (AM), on the other hand, is on the simulated (Monte Carlo based)
environment.
It is used in this paper for testing purposes.
In real world applications, the observation process is the
actual process obtained from real physical process.
To test the model robustness, we consider the case when the NM's observation
noise differs from the AM's observation noise.
In particular, we consider the following two models:
\[
\left\{\begin{array}{l}
{(\NM)}:\
  \left\{\begin{array}{l}
x_{n+1}=(1+0.1\eta)x_n+\sqrt{\eta}\ \sigma u_n,\ x_0=x,\\
y_n= x_n+\sigma_0^{\NM}v_n, \ n=0,1,2,\ldots,
\end{array}\right.\\
  {(\AM)}:\
  \left\{\begin{array}{l}
x_{n+1}=(1+0.1\eta)x_n+\sqrt{\eta}\ \sigma u_n,\ x_0=x,\\
y_n= x_n+\sigma_0^{\AM}v_n, \ n=0,1,2,\ldots,
\end{array}\right.
\end{array}\right.
\]

First, with fixed $\sigma_0^\AM=0.5$, we vary $\sigma_0^\NM$.
The corresponding DF errors and the KF errors are given in
Table~\ref{sigma-NM}. The KF depends heavily on
nominal observation noise. On the other hand, the DF is more robust when $\sigma_0^\NM\geq 0.5$.
Also, the DF needs the nominal observation noise to
be in normal range (not too small) in order to properly train the DNN.
Namely, some noise is necessary when training a DNN.
Or a noise process in fact helps in the training stage of a DNN.

\begin{table}
\begin{center}
    {\small
  \begin{tabular}{|c|c|c|c|c|c|c|}\hline
$\sigma_0^\NM$ & 0.1 & 0.5 & 1.0 & 1.5 & 2.0 & 2.5 \\ \hline
DF &  8.89 & 5.71 & 5.64 & 5.62 & 5.62 & 5.62   \\ \hline
KF &  6.54 & 4.08 & 4.59 & 5.33 & 6.04 & 6.70   \\ \hline
  \end{tabular}
  }
\medskip
\caption{Error dependence on $\sigma_0^\NM$.}\label{sigma-NM}
    \end{center}
\end{table}

Next, with fixed $\sigma_0^\NM=0.5$, we vary $\sigma_0^\AM$.
As the actual observation noise increases, both DF and KF deteriorate and the corresponding
errors increase as shown
in Table~\ref{sigma-AM}. The DF appears to be more robust than the KF because
it is less sensitive to such changes than the KF.

\begin{table}
  \begin{center}
        {\small
\begin{tabular}{|c|c|c|c|c|c|c|}\hline
$\sigma_0^\AM$ & 0.1 & 0.5 & 1.0 & 1.5 & 2.0 & 2.5 \\ \hline
DF &  5.25 & 5.71 & 6.99 & 8.71 & 10.63 & 12.65   \\ \hline
KF &  2.90 & 4.08 & 6.49 & 9.15 & 11.88 & 14.59   \\ \hline
\end{tabular}
}
\medskip
\caption{Error dependence on $\sigma_0^\NM$.}\label{sigma-AM}
    \end{center}
\end{table}

\subsection{Nonlinear Models}
In this section, we consider nonlinear (NL) models and comparison of the DF with
the corresponding extended Kalman filter.
We consider the two (NM and AM) models:
\[
\left\{\begin{array}{l}
{(\NM)}:\
  \left\{\begin{array}{l}
x_{n+1}=x_n+\eta \sin(5x_n)+\sqrt{\eta}\ \sigma u_n,\ x_0=x,\\
y_n= x_n+\sigma_0^{\NM}v_n, \ n=0,1,2,\ldots,
\end{array}\right.\\
  {(\AM)}:\
  \left\{\begin{array}{l}
  x_{n+1}=x_n+\eta \sin(5x_n)+\sqrt{\eta}\ \sigma u_n,\ x_0=x,\\
  y_n= x_n+\sigma_0^{\AM}v_n, \ n=0,1,2,\ldots,
\end{array}\right.
\end{array}\right.
\]

We take $T=5$, step size $\eta=0.005$, and $N=1000$. We also take
$\sigma=0.7$ and $\sigma_0^\AM=0.5$.
With these specifications, we vary $\sigma_0^\NM$.
In Table~\ref{N2N-NM}, it can be seen that deep filter is more robust and
less dependent on nominal observation noise changes when compared against
the corresponding extended Kalman filter.
We also note that when training the DF, the observation noise in
training data should not be too small. This is typical in DNN training.
Too little noise will not provide necessary variations when training the DNN.

\begin{table}
  \begin{center}
        {\small
\begin{tabular}{|c|c|c|c|c|c|c|}\hline
$\sigma_0^\NM$ & 0.1 & 0.5 & 1.0 & 1.5  & 2.0 & 2.5 \\ \hline
DF &  12.24 & 7.75 & 7.62 & 7.58 & 7.56  & 7.56   \\ \hline
EKF &  9.22  & 5.58 & 6.63 & 8.29 & 10.13 & 12.14   \\ \hline
\end{tabular}
}
\medskip
\caption{(NM=NL, AM=NL): Error dependence on $\sigma_0^\NM$.}\label{N2N-NM}
    \end{center}
\end{table}

Next, we fix $\sigma_0^\NM=0.5$ and vary $\sigma_0^\AM$,
Increasing in actual observation noise will make filtering more
difficult and increase the corresponding filtering errors.
This is confirmed in Table~\ref{N2N-AM}. Also, the DF is less
affected than the EKF as $\sigma_0^\AM$ increases.

\begin{table}
  \begin{center}
        {\small
\begin{tabular}{|c|c|c|c|c|c|c|}\hline
$\sigma_0^\AM$ & 0.1 & 0.5 & 1.0 & 1.5  & 2.0 & 2.5 \\ \hline
DF &  7.01  & 7.75 & 9.68 &12.22 & 15.01 & 17.89   \\ \hline
EKF &  4.19  & 5.58 & 8.62 &12.26 & 16.17 & 20.07   \\ \hline
\end{tabular}
}
\medskip
\caption{(NM=NL, AM=NL): Error dependence on $\sigma_0^\AM$.}\label{N2N-AM}
    \end{center}
\end{table}

\subsubsection*{Mixed Nonlinear and Linear Models}
In this section, we consider the case when the linearity of the NM and the AM differs.
First, we use a linear model to train the DNN while the actual model
is in fact nonlinear. In particular, we consider the following two models:
\[
\left\{\begin{array}{l}
{(\NM)}:\
  \left\{\begin{array}{l}
x_{n+1}=(1+0.1\eta)x_n+\sqrt{\eta}\ \sigma u_n,\ x_0=x,\\
y_n= x_n+\sigma_0^{\NM}v_n, \ n=0,1,2,\ldots,
\end{array}\right.\\
  {(\AM)}:\
  \left\{\begin{array}{l}
  x_{n+1}=x_n+\eta \sin(5x_n)+\sqrt{\eta}\ \sigma u_n,\ x_0=x,\\
  y_n= x_n+\sigma_0^{\AM}v_n, \ n=0,1,2,\ldots,
\end{array}\right.
\end{array}\right.
\]

  In this case, we fix $\sigma_0^\AM=0.5$ and vary $\sigma_0^\NM$,
  The DF is barely affected with the training model (NM) when $\sigma_0^\NM\geq0.5$ as shown
  in Table~\ref{L2N-NM}.
  The dependence of the KF errors on $\sigma_0^\NM$ is more pronounced though.

  Then we fix $\sigma_0^\NM=0.5$ and vary $\sigma_0^\AM$.
  It can be seen from
  Table~\ref{L2N-AM} that both
  the KF and DF errors increase in $\sigma_0^\AM$, but
  the DF is more robust because its errors are less sensitive.
\begin{table}
  \begin{center}
        {\small
\begin{tabular}{|c|c|c|c|c|c|c|}\hline
$\sigma_0^\NM$ & 0.1 & 0.5 & 1.0 & 1.5  & 2.0 & 2.5 \\ \hline
DF &  12.91 & 7.75 & 7.61 & 7.57 & 7.56  & 7.56   \\ \hline
KF &  9.30  & 5.70 & 6.22 & 7.05 & 7.82  & 8.58   \\ \hline
\end{tabular}
}
\medskip
\caption{(NM=L, AM=NL): Error dependence on $\sigma_0^\NM$.}\label{L2N-NM}
    \end{center}
\end{table}

\begin{table}
\begin{center}
    {\small
  \begin{tabular}{|c|c|c|c|c|c|c|}\hline
$\sigma_0^\AM$ & 0.1 & 0.5 & 1.0 & 1.5  & 2.0 & 2.5 \\ \hline
DF &  7.00  & 7.75 & 9.72 &12.27 & 15.05 & 17.88   \\ \hline
KF &  3.98  & 5.70 & 9.16 &12.96 & 16.77 & 20.48   \\ \hline
\end{tabular}
    }
\medskip
    \caption{(NM=L, AM=NL): Error dependence on $\sigma_0^\AM$.}\label{L2N-AM}
    \end{center}
\end{table}

Finally, we consider the case when the training model is nonlinear while
the actual model is linear. We consider the following two models:
\[
\left\{\begin{array}{l}
  {(\NM)}:\
  \left\{\begin{array}{l}
  x_{n+1}=x_n+\eta \sin(5x_n)+\sqrt{\eta}\ \sigma u_n,\ x_0=x,\\
y_n= x_n+\sigma_0^{\NM}v_n, \ n=0,1,2,\ldots,
\end{array}\right.\\
  {(\AM)}:\
  \left\{\begin{array}{l}
x_{n+1}=(1+0.1\eta)x_n+\sqrt{\eta}\ \sigma u_n,\ x_0=x,\\
  y_n= x_n+\sigma_0^{\AM}v_n, \ n=0,1,2,\ldots,
  \end{array}\right.
  \end{array}\right.
  \]

  Such changes help to improve the filtering outcomes in both DF and EKF
  when we vary $\sigma_0^\NM$ and  $\sigma_0^\AM$. It appears
   that this helps the DF
  more in reduction of estimation errors as shown in Tables~\ref{N2L-NM} and
  \ref{N2L-AM} when compared with Tables~\ref{L2N-NM} and \ref{L2N-AM}, respectively.
\begin{table}
  \begin{center}
        {\small
\begin{tabular}{|c|c|c|c|c|c|c|}\hline
$\sigma_0^\NM$ & 0.1 & 0.5 & 1.0 & 1.5   & 2.0 & 2.5 \\ \hline
DF &  8.53  & 5.75 & 5.66 & 5.64 & 5.63  & 5.62   \\ \hline
KF &  6.54  & 4.24 & 5.34 & 7.09 & 9.27  & 11.94  \\ \hline
\end{tabular}
}
\medskip
\caption{(NM=NL, AM=L): Error dependence on $\sigma_0^\NM$.}\label{N2L-NM}
    \end{center}
\end{table}

\begin{table}
  \begin{center}
        {\small
\begin{tabular}{|c|c|c|c|c|c|c|}\hline
$\sigma_0^\AM$ & 0.1 & 0.5 & 1.0 & 1.5  & 2.0 & 2.5 \\ \hline
DF &  5.28  & 5.75 & 7.02 &8.72 & 10.64 & 12.66   \\ \hline
KF &  3.12  & 4.24 & 6.55 &9.13 & 11.80 & 14.50   \\ \hline
\end{tabular}
}
\medskip
\caption{(NM=NL, AM=L): Error dependence on $\sigma_0^\AM$.}\label{N2L-AM}
    \end{center}
\end{table}

\subsection{A Randomly Switching Model}
Finally, we consider a switching model with jumps and apply the DF
to these models. Note that
neither the KF nor the EKF can be used
in this case due to the presence of the switching process.

As a demonstration,
we consider the following model:
\[
\left\{\begin{array}{l}
(\NM): \left\{\begin{array}{l}
x_n=\sin(n\eta\al_n+\sigma u_n),\ n=0,1,2,\ldots,\\
y_n=x_n+\sigma_0^\NM v_n,
\end{array}\right.\\
(\AM):\left\{\begin{array}{l}
x_n=\sin(n\eta\al_n+\sigma u_n),\ n=0,1,2,\ldots.\\
y_n=x_n+\sigma_0^\AM v_n.
\end{array}\right.
\end{array}\right.
\]

\begin{center}
\begin{figure}[h!tb]
\begin{center}
\mbox{\subfigure[State $x_n$ (sample path 1)]{
         \psfig{figure=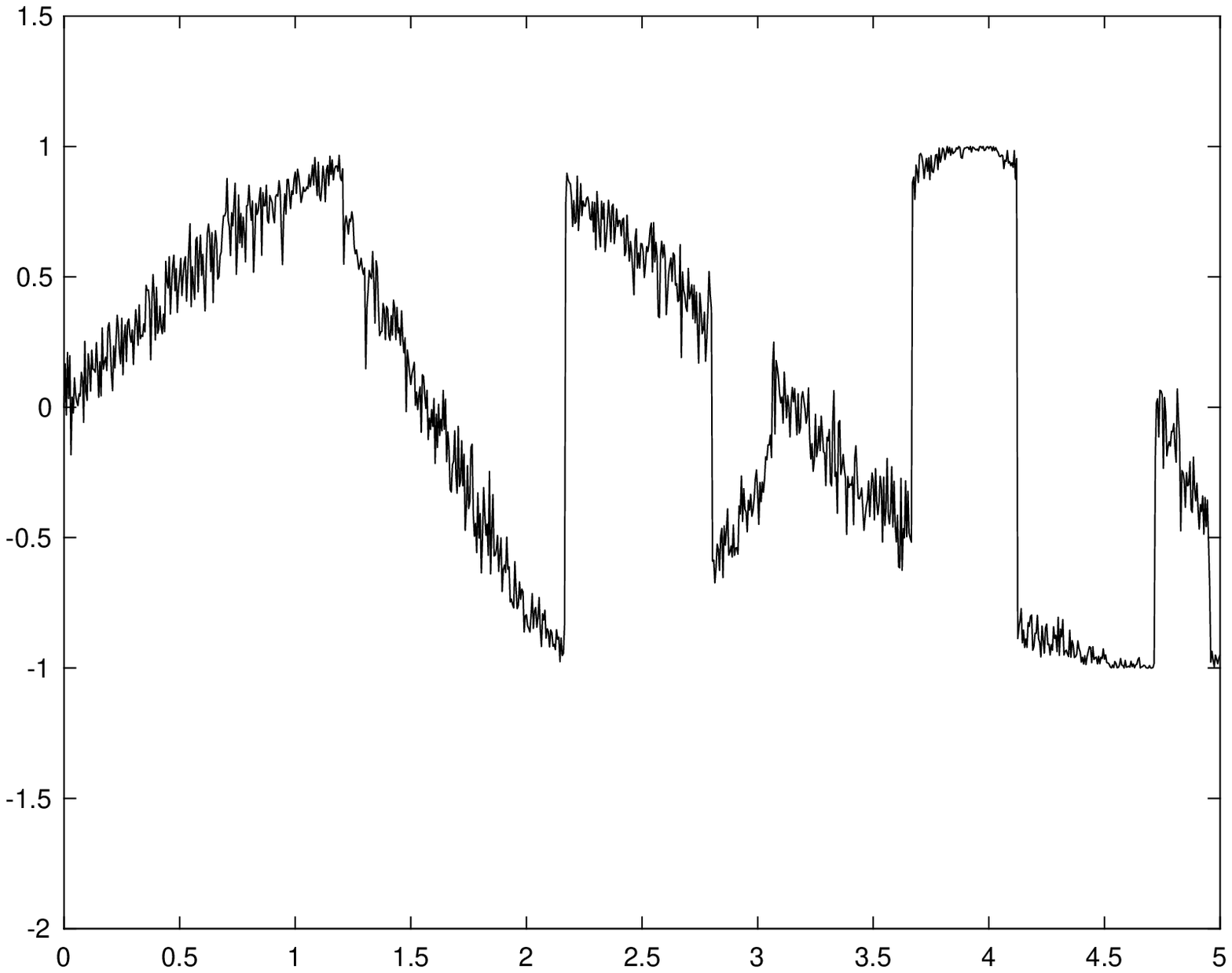,width=0.35\linewidth }} \qquad
      \subfigure[State $x_n$ (sample path 2)]{
        \psfig{figure=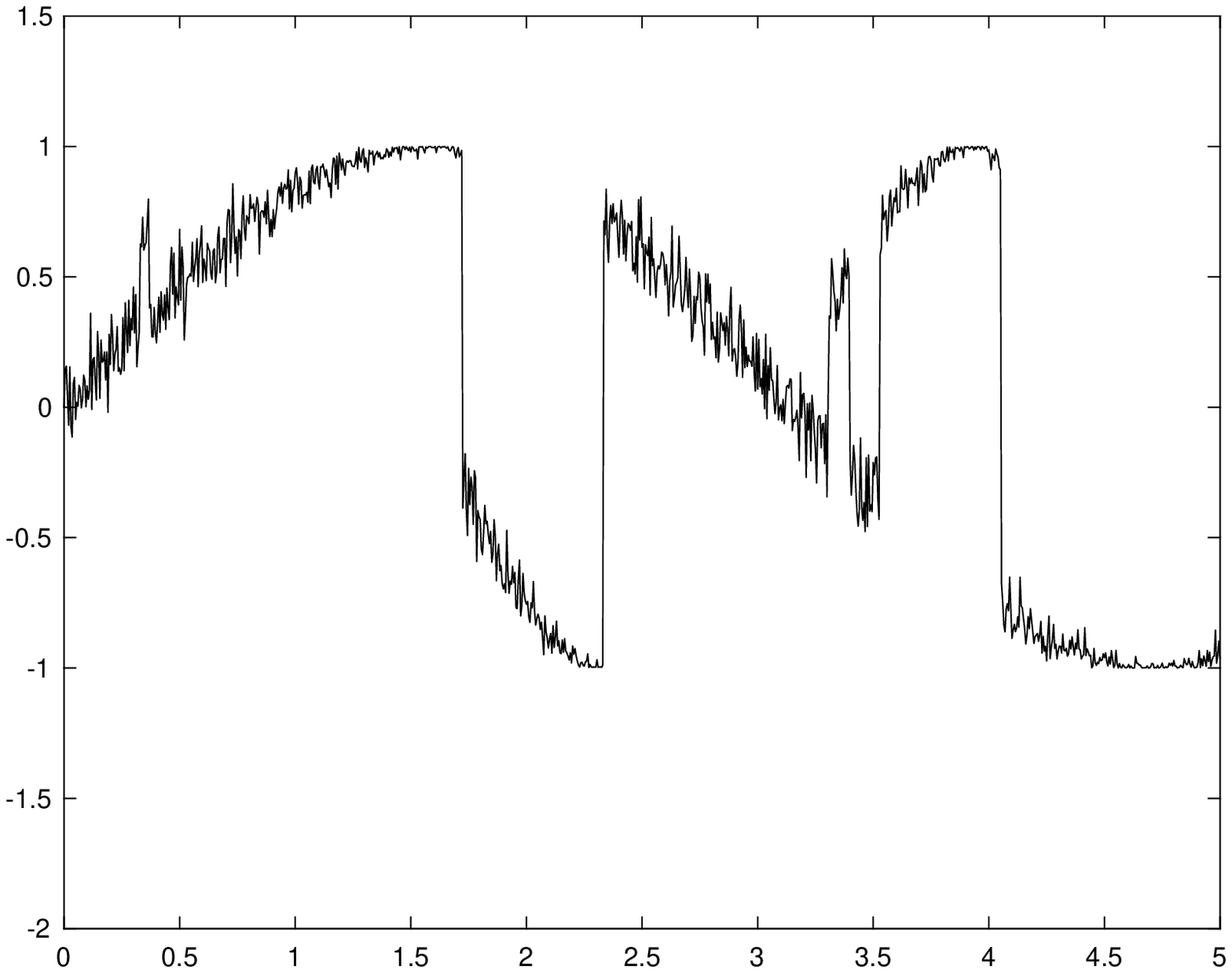,width=0.35\linewidth }}}

\mbox{\subfigure[DF $\wdt x_n$ (sample path 1)]{
    \psfig{figure=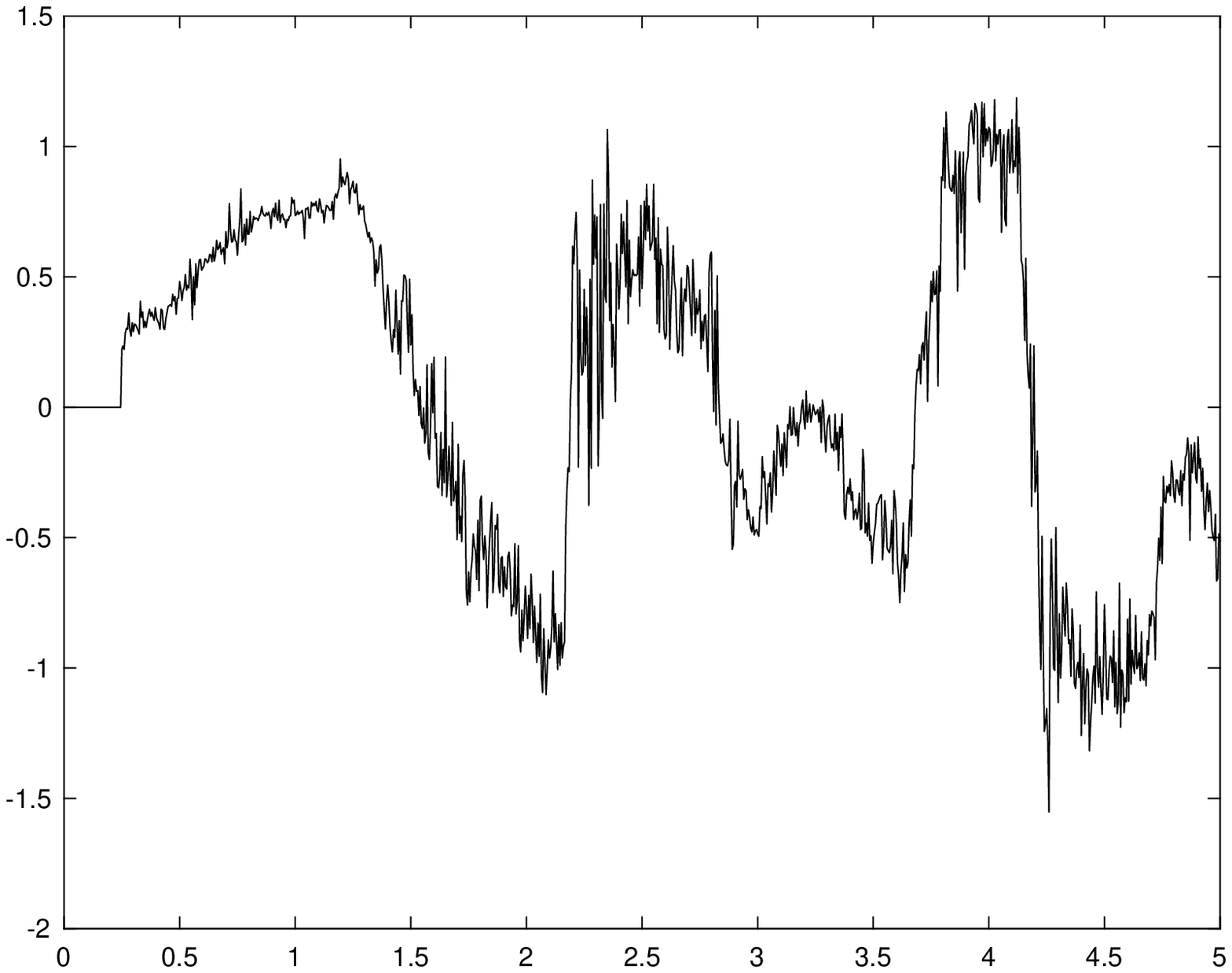,width=0.35\linewidth }} \qquad
      \subfigure[DF $\wdt x_n$ (sample path 2)]{
        \psfig{figure=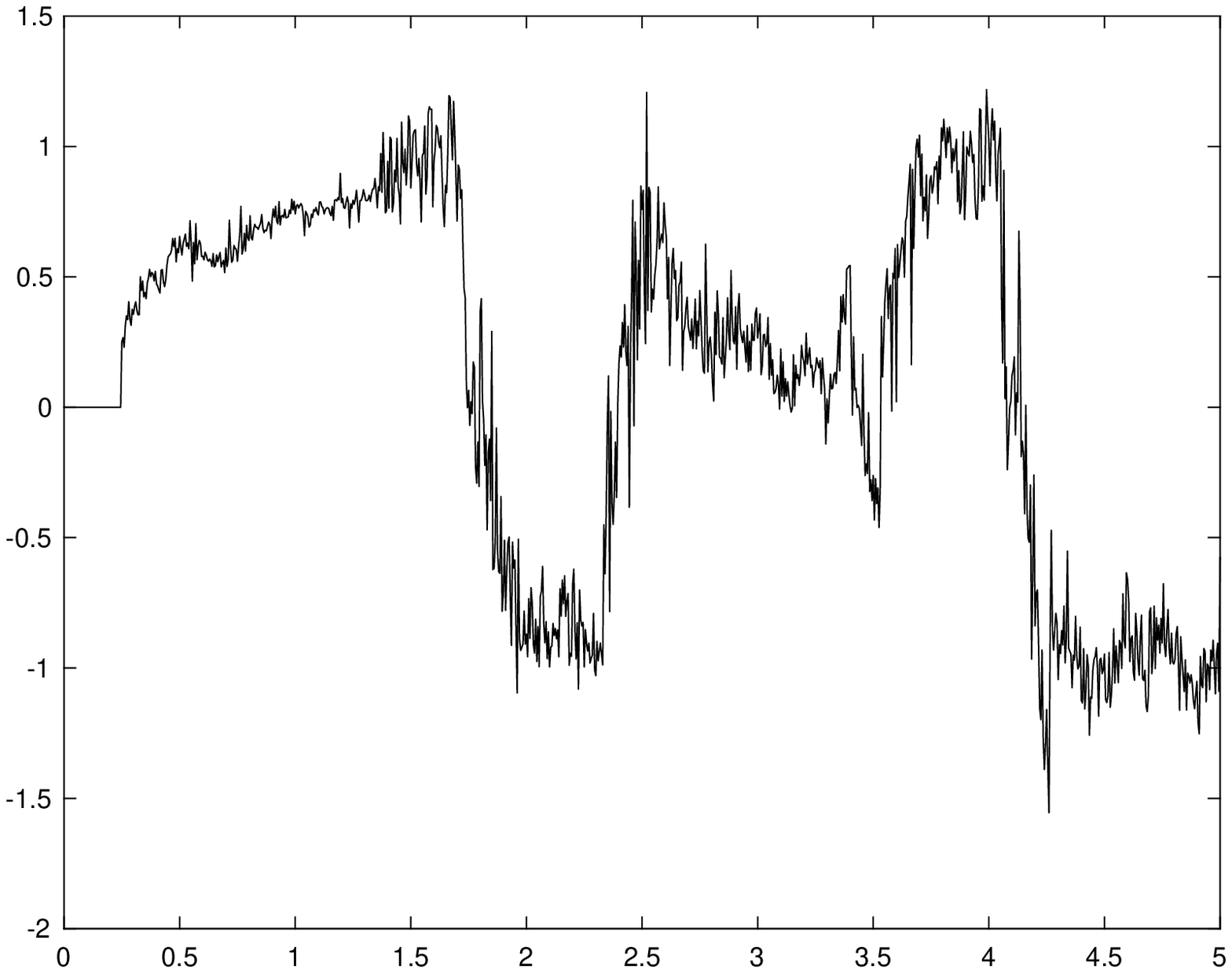,width=0.35\linewidth }}}

\mbox{\subfigure[Error $x_n-\wdt x_n$ (sample path 1)]{
         \psfig{figure=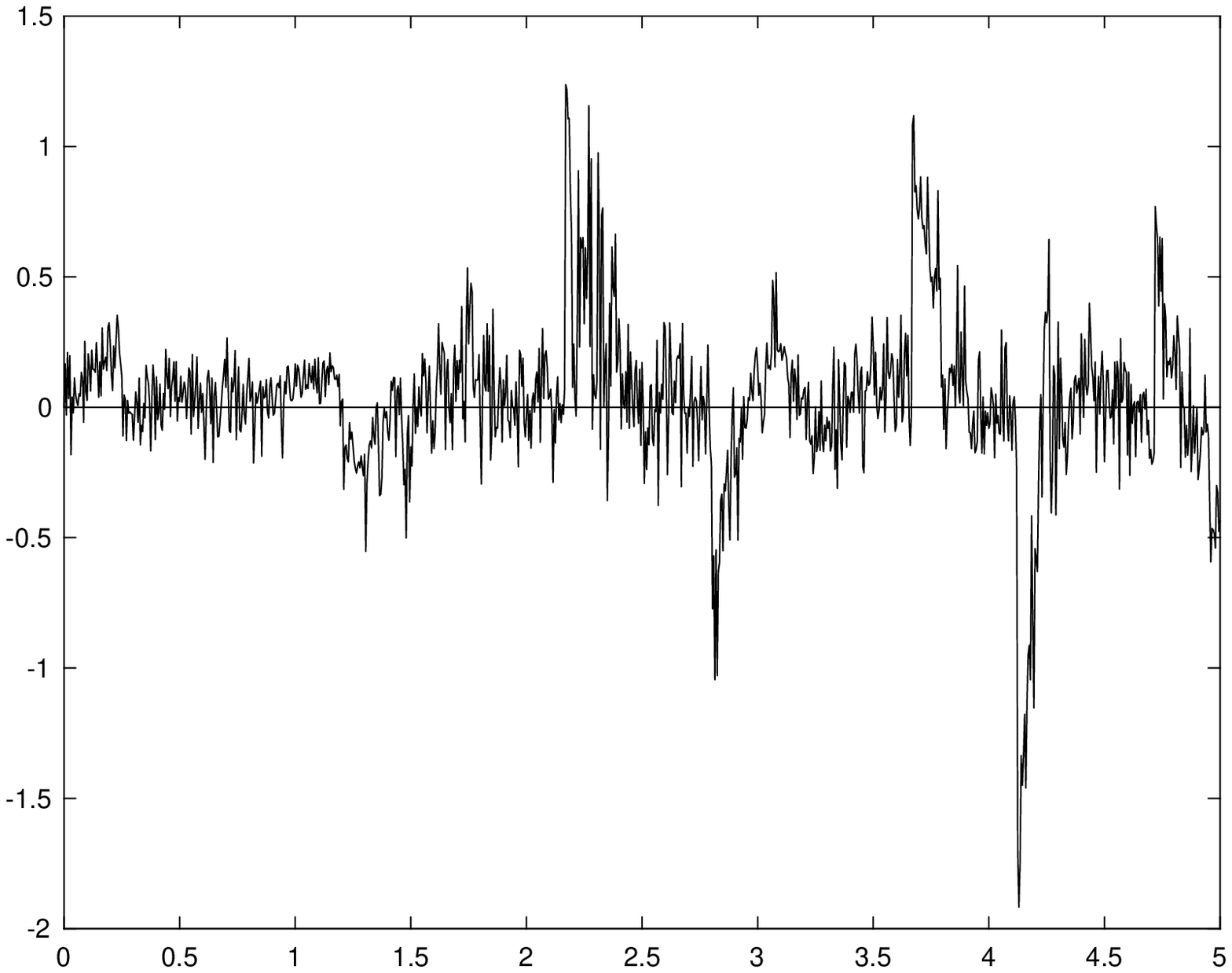,width=0.35\linewidth }} \qquad
      \subfigure[Error $x_n-\wdt x_n$ (sample path 2)]{
         \psfig{figure=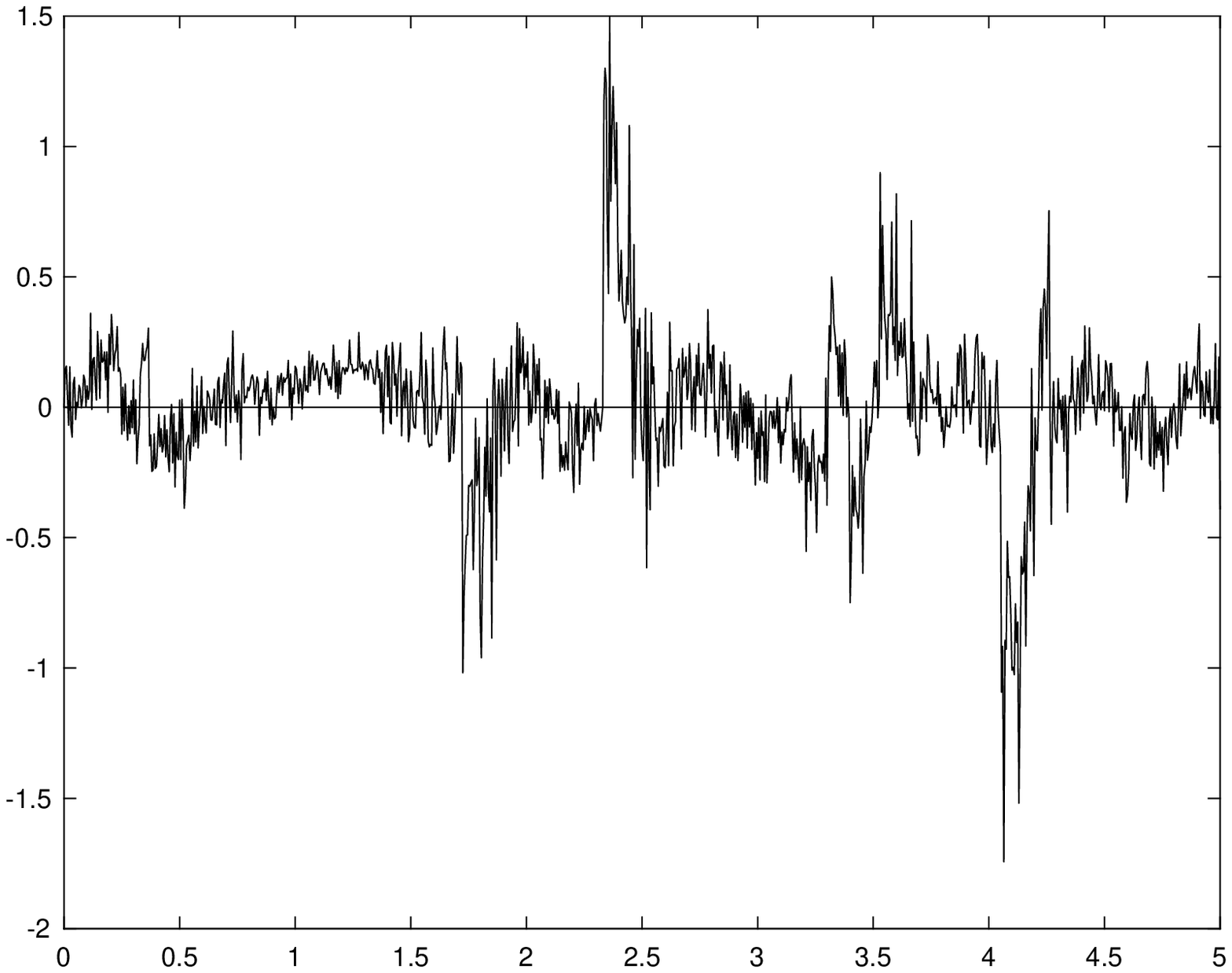,width=0.35\linewidth }}}
\end{center}
\medskip
\caption{\footnotesize Deep filter for randomly switched system:
Two sample paths of $x_n$, $\wdt x_n$, and errors.}\label{switching}
  \end{figure}
\end{center}

We take $\al(t)\in\{1,2\}$ to be a continuous-time Markov chain
with generator $Q=\(\begin{array}{cc}
-2 & 2\\
2 & -2\\
\end{array}\)$;
see, for example, Yin and Zhang \cite{YinZbook} for details.
Using step size $\eta=0.005$ to discretize $\al(t)$ to get
$\al_n=\al(n\eta)$.
We also take $\sigma=0.1$ and $\sigma_0=0.3$.
Two sample paths of $x_n$, $\wdt x_n$, and the corresponding errors are plotted in
Figure~\ref{switching}.

The DF appears to be effective and it catches up quickly the jumps of $x_n$.
Then, in Tables~\ref{Switching-NM} and \ref{Switching-AM}, we provide
the errors when one of $\sigma_0^\NM$ and $\sigma_0^\AM$ is fixed to 0.3 and the
other varies. These errors are larger
than that
of linear and nonlinear models
in the previous section mainly because of the presence of jumps.
In addition, as $\sigma_0^\NM$ moves away from 0.2-0.3, the errors increase.
Similarly as in the previous linear and nonlinear models, the errors
increase in $\sigma_0^\AM$ when $\sigma_0^\NM$ is fixed.
Overall, the DF shows strong adaptiveness and effectiveness
in filtering under highly nonlinear with switching (jumps) dynamic models.
\begin{table}
  \begin{center}
        {\small
\begin{tabular}{|c|c|c|c|c|c|c|}\hline
$\sigma_0^\NM$ & 0.1 & 0.2 & 0.3 & 0.4  & 0.5 & 0.6 \\ \hline
DF &  14.88  & 13.41 & 13.78 & 15.05 & 16.33 & 17.72   \\ \hline
\end{tabular}
}
\medskip
\caption{Switching model: Error dependence on $\sigma_0^\NM$.}\label{Switching-NM}
    \end{center}
\end{table}

\begin{table}
\begin{center}
    {\small
  \begin{tabular}{|c|c|c|c|c|c|c|}\hline
$\sigma_0^\AM$ & 0.1 & 0.2 & 0.3 & 0.4  & 0.5 & 0.6 \\ \hline
DF &  12.22  & 12.74 & 13.78 & 15.19 & 16.82 & 18.59   \\ \hline
  \end{tabular}
  }
\medskip
\caption{Switching model: Error dependence on $\sigma_0^\AM$.}\label{Switching-AM}
    \end{center}
\end{table}

\section{Concluding Remarks}\label{sec:con}
In this paper, we developed a new approach using
deep learning for
stochastic filtering.
We explore
deep neural networks by providing preliminary experiments on
various dynamic models.
Naturally one would be interested in any theoretical analysis on related convergence,
extensive numerical tests on high-dimensional models with possible high nonlinearity,
any genuine real world applications. All in all,
this paper raises some opportunities and challenges. Nevertheless, there are more questions than answers.


\end{document}